\input amstex
\input amsppt.sty
\magnification=\magstep1
\hsize=30truecc
\vsize=22.2truecm
\baselineskip=16truept
\NoBlackBoxes
\TagsOnRight \pageno=1 \nologo
\def\Z{\Bbb Z}
\def\N{\Bbb N}

\def\l{\left}
\def\r{\right}
\def\bg{\bigg}
\def\({\bg(}
\def\[{\bg\lfloor}
\def\){\bg)}
\def\]{\bg\rfloor}
\def\t{\text}
\def\f{\frac}

\def\bi{\binom}
\def\eq{\equiv}

\def\ls{\leqslant}
\def\gs{\geqslant}
\def\mo{\roman{mod}}

\def\al{\alpha}
\def\da{\delta}

\def\Proof{\noindent{\it Proof}}

\def\Remark{\medskip\noindent{\it  Remark}}

\def\Ack{\medskip\noindent {\bf Acknowledgments}}
\hbox {Publ. Math. Debrecen 79(2011), in press.}
\bigskip
\topmatter
\title On harmonic numbers and Lucas sequences\endtitle
\author Zhi-Wei Sun\endauthor
\leftheadtext{Zhi-Wei Sun} \rightheadtext{On harmonic numbers and Lucas sequences}
\affil Department of Mathematics, Nanjing University\\
 Nanjing 210093, People's Republic of China
  \\  zwsun\@nju.edu.cn
  \\ {\tt http://math.nju.edu.cn/$\sim$zwsun}
\endaffil
\abstract Harmonic numbers $H_k=\sum_{0<j\ls k}1/j\
(k=0,1,2,\ldots)$ arise naturally in many fields of mathematics. In
this paper we initiate the study of congruences involving both harmonic
numbers and Lucas sequences. One of our three theorems is as follows: Let
$u_0=0,\ u_1=1,$ and $u_{n+1}=u_n-4u_{n-1}$ for $n=1,2,3,\ldots$.
Then, for any prime $p>5$ we have
$$\sum_{k=0}^{p-1}\f {H_k}{2^k}u_{k+\da}\eq0\ (\mo\ p),$$
where $\da=0$ if $p\eq1,2,4,8\ (\mo\ 15)$, and $\da=1$ otherwise.
\endabstract
\thanks 2010 {\it Mathematics Subject Classification}.\,Primary 11A07, 11B39;
Secondary 05A10, 33B99.
\newline\indent Supported by the National Natural Science
Foundation (grant 10871087) and the Overseas Cooperation Fund (grant 10928101) of China.
\endthanks
\endtopmatter
\document

\heading{1. Introduction}\endheading

Harmonic numbers are those rational numbers given by
$$H_n=\sum_{0<k\ls n}\f1k\quad \ (n\in\N=\{0,1,2,\ldots\}).$$
They play important roles in mathematics; see, e.g., [BPQ] and [BB].

In 1862 J. Wolstenholme [W] (see also [HT]) discovered that for any prime $p>3$ we have
$$H_{p-1}=\sum_{k=1}^{p-1}\f1k\eq0\ (\mo\ p^2).$$
In a previous paper [Su] the
author developed the arithmetic theory of harmonic numbers by
proving the following fundamental congruences for primes $p>3$:
$$\sum_{k=1}^{p-1}H_k^2\eq 2p-2\ (\mo\ p^2),
\ \ \ \sum_{k=1}^{p-1}H_k^3\eq 6\ (\mo\ p),$$ and
$$\sum_{k=1}^{p-1}\f{H_k^2}{k^2}\eq0\ (\mo\ p)\qquad \t{provided}\ p>5.$$

In this paper we initiate the investigation of congruences involving both harmonic numbers and
Lucas sequences.

For $A,B\in\Z$ the Lucas sequences $u_n=u_n(A,B)\ (n\in\N)$ and
$v_n=v_n(A,B)\ (n\in\N)$ are defined as follows:
$$u_0=0,\ u_1=1,\ \t{and}\ u_{n+1}=Au_n-Bu_{n-1}\ (n=1,2,3,\ldots);$$
$$v_0=2,\ v_1=A,\ \t{and}\ v_{n+1}=Av_n-Bv_{n-1}\ (n=1,2,3,\ldots).$$
The sequence $\{v_n\}_{n\gs0}$ is called the companion sequence of
$\{u_n\}_{n\gs0}$. The characteristic equation $x^2-Ax+B=0$ of the
sequences $\{u_n\}_{n\gs0}$ and $\{v_n\}_{n\gs0}$ has two roots
$$\al=\f{A+\sqrt{\Delta}}2\quad\t{and}\quad\beta=\f{A-\sqrt{\Delta}}2,$$
where $\Delta=A^2-4B$. It is well known that for any $n\in\N$ we
have
$$Au_n+v_n=2u_{n+1},\ (\al-\beta)u_n=\al^n-\beta^n,\ \ \t{and}\ \  v_n=\al^n+\beta^n.$$
(See, e.g., [R, pp.\,41--44].)
Note that those $F_n=u_n(1,-1)$ and $L_n=v_n(1,-1)$ are well-known
Fibonacci numbers and Lucas numbers respectively.

 Here is our first theorem.

\proclaim{Theorem 1.1} Let $p>3$ be a prime and let $A,B\in\Z$ with $p\nmid A$.
Then
$$\sum_{k=1}^{p-1}\f{v_k(A,B)H_k}{kA^k}\eq0\ (\mo\ p)\tag1.1$$
and
$$\sum_{k=1}^{p-1}\f{u_k(A,B)H_k}{kA^k}\eq\f 2p\sum_{k=1}^{p-1}\f{u_k(A,B)}{kA^k}\ (\mo\ p).\tag1.2$$
\endproclaim

Since $v_k(2,1)=2$ for all $k\in\N$, Theorem 1.1 yields the following consequence.

\proclaim{Corollary 1.1 {\rm ([Su])}} For any prime $p>3$ we have
$$\sum_{k=1}^{p-1}\f{H_k}{k2^k}\eq0\ (\mo\ p).$$
\endproclaim
\Remark. In 1987 S. W. Coffman [C] proved that $\sum_{k=1}^\infty H_k/(k2^k)=\pi^2/12$.

\medskip

Applying Theorem 1.1 with $A=1$ and $B=-1$ we get the following corollary.
\proclaim{Corollary 1.2} Let $p>3$ be a prime. Then
$$\sum_{k=1}^{p-1}\f{H_kL_k}k\eq0\ (\mo\ p)\ \ \t{and}\ \ \sum_{k=1}^{p-1}\f{F_kH_k}k\eq\f 2p\sum_{k=1}^{p-1}\f{F_k}k
\ (\mo\ p).\tag1.3$$
\endproclaim

Let $\omega$ denote the cubic root $(-1+\sqrt{-3})/2$. For $n\in\N$ we have
$$u_n(-1,1)=u_n(\omega+\bar\omega,\omega\bar\omega)=\f{\omega^n-\bar\omega^n}{\sqrt{-3}}=\l(\f n3\r)$$
and
$$u_n(1,1)=(-1)^{n-1}u_n(-1,1)=(-1)^{n-1}\l(\f n3\r),$$
where $(-)$ denotes the Jacobi symbol. By induction, for any $k\in\N$ we have
$$u_{4k}(2,2)=0,\ u_{4k+1}(2,2)=(-4)^k\ \t{and}\ u_{4k+2}(2,2)=u_{4k+3}(2,2)=2(-4)^k.$$

Now we state our second theorem.

\proclaim{Theorem 1.2} Let $p>3$ be a prime.

{\rm (i)} Let $A,B\in\Z$ with $p\nmid B$ and $(\f{A^2-4B}p)=1$. Then, for any $n\in\N$ we have
$$\sum_{k=0}^{p-1}(1+B^{-k})u_k(A,B)H_k^n\eq0\ (\mo\ p)\tag1.4$$
and
$$ \sum_{k=0}^{p-1}(1-B^{-k})v_k(A,B)H_k^n\eq0\ (\mo\ p).\tag1.5$$

{\rm (ii)} We have
$$\sum_{k=0}^{p-1}(-1)^k\l(\f k3\r)H_k\eq0\ (\mo\ p)\tag1.6$$
and $$\sum_{k=0}^{p-1}\l(\f k3\r)H_k\eq\f{(\f p3)-1}4q_p(3)\ (\mo\ p),\tag1.7$$
where $q_p(3)$ refers to the Fermat quotient $(3^{p-1}-1)/p$.
Also,
$$\sum_{k=0}^{p-1}(-1)^k\l(\f k3\r)kH_k\eq\f{1-(\f p3)}2\ (\mo\ p)\tag1.8$$
and
$$\sum_{k=0}^{p-1}(1+2^{-k})u_k(2,2)H_k\eq0\ (\mo\ p).\tag1.9$$
\endproclaim

Since $F_{2k}=u_k(3,1)$, $F_k=u_k(1,-1)$ and $L_k=v_k(1,-1)$ for all $k\in\N$, Theorem 1.2(i)
implies the following result.

\proclaim{Corollary 1.3} Let $p$ be a prime with $p\eq\pm1\ (\mo\ 5)$. Then
$$\sum_{k=0}^{p-1}F_{2k}H_k^n\eq\sum^{p-1}\Sb k=0\\2\mid k\endSb F_kH_k^n\eq\sum^{p-1}\Sb k=0\\2\nmid k\endSb L_kH_k^n
\eq0\ (\mo\ p)\tag1.10$$
for every $n=0,1,2,\ldots$.
\endproclaim

Our third theorem seems curious and unexpected.

 \proclaim{Theorem 1.3} Let $p>5$ be a prime.
If $(\f p{15})=1$, i.e., $p\eq 1,2,4,8\ (\mo\ 15)$, then
$$\sum_{k=0}^{p-1}\f{u_k(1,4)}{2^k}H_k\eq0\ (\mo\ p).\tag1.11$$
If $(\f p{15})=-1$, i.e., $p\eq 7,11,13,14\ (\mo\ 15)$, then
$$\sum_{k=0}^{p-1}\f{u_{k+1}(1,4)}{2^k}H_k\eq0\ (\mo\ p).\tag1.12$$
\endproclaim

In the next section we are going to prove Theorems 1.1 and 1.2.
Section 3 is devoted to the proof of Theorem 1.3.

To conclude this section we pose three related conjectures.

Recall that harmonic numbers of the second order
are given by
$$H_n^{(2)}=\sum_{0<k\ls n}\f1{k^2}\quad \ (n=0,1,2,\ldots).$$

\proclaim{Conjecture 1.1} Let $p>3$ be a prime. Then
$$\sum_{k=0}^{p-1}(-2)^k\bi{2k}kH_k^{(2)}\eq\f23q_p(2)^2\ (\mo\ p)$$
where $q_p(2)=(2^{p-1}-1)/p$. If $p>5$, then we have
$$\sum_{k=0}^{p-1}(-1)^k\bi{2k}kH_k^{(2)}\eq\f 52\l(\f p5\r)\f {F^2_{p-(\f p5)}}{p^2}\ (\mo\ p).$$
\endproclaim

\proclaim{Conjecture 1.2} Let $p$ be an odd prime. Then
$$\align\sum_{k=0}^{p-1}\f{u_k(2,-1)}{(-8)^k}\bi{2k}k^2\eq&0\ (\mo\ p)\ \ \ \t{if}\ p\eq 5\ (\mo\ 8),
\\\sum_{k=0}^{p-1}\f{u_k(2,-1)}{32^k}\bi{2k}k^2\eq&0\ (\mo\ p)\ \ \ \t{if}\ p\eq 7\ (\mo\ 8),
\\\sum_{k=0}^{p-1}\f{v_k(2,-1)}{(-8)^k}\bi{2k}k^2\eq&0\ (\mo\ p)\ \ \ \t{if}\ p\eq5,7\ (\mo\ 8),
\\\sum_{k=0}^{p-1}\f{v_k(2,-1)}{32^k}\bi{2k}k^2\eq&0\ (\mo\ p)\ \ \ \t{if}\ p\eq 5\ (\mo\ 8).
\endalign$$
Also,
$$\gather\sum_{k=0}^{p-1}\f{u_k(4,1)}{4^k}\bi{2k}k^2\eq0\ (\mo\ p)\ \ \ \t{if}\ p\eq 2\ (\mo\ 3),
\\\sum_{k=0}^{p-1}\f{u_k(4,1)}{64^k}\bi{2k}k^2\eq0\ (\mo\ p)\ \ \ \t{if}\ p\eq 11\ (\mo\ 12),
\\\sum_{k=0}^{p-1}\f{v_k(4,1)}{4^k}\bi{2k}k^2\eq\sum_{k=0}^{p-1}\f{v_k(4,1)}{64^k}\bi{2k}k^2\eq0\ (\mo\ p)
 \ \t{if}\ p\eq 5\ (\mo\ 12).
 \endgather$$
\endproclaim

\proclaim{Conjecture 1.3} Let $p>3$ be a prime. Then
$$\sum_{k=0}^{p-1}\bi{p-1}k\bi{2k}k((-1)^k-(-3)^{-k})\eq\l(\f p3\r)(3^{p-1}-1)\ (\mo\ p^3).$$
If $p\eq\pm1\ (\mo\ 12)$, then
$$\sum_{k=0}^{p-1}\bi{p-1}k\bi{2k}k(-1)^ku_k(4,1)\eq(-1)^{(p-1)/2}u_{p-1}(4,1)\ (\mo\ p^3).$$
If $p\eq\pm1\ (\mo\ 8)$, then
$$\sum_{k=0}^{p-1}\bi{p-1}k\bi{2k}k\f{u_k(4,2)}{(-2)^k}\eq(-1)^{(p-1)/2}u_{p-1}(4,2)\ (\mo\ p^3).$$
\endproclaim

\heading{2. Proofs of Theorems 1.1 and 1.2}\endheading

Let $p$ be an odd prime. For $k=0,1,\ldots,p-1$ we obviously have
$$(-1)^k\bi{p-1}k=\prod_{0<j\ls k}\l(1-\f pj\r)\eq1-pH_k\ (\mo\ p^2).\tag2.1$$
This basic fact is useful in the study of congruences involving harmonic numbers.

\proclaim{Lemma 2.1} Let $n\gs j>0$ be integers. Then
$$\sum_{k=j}^n\bi{k-1}{j-1}=\bi{n}j.$$
\endproclaim
\Proof. This is a known identity due to Shih-chieh Chu (cf. (5.26) of [GKP, p.\,169]). By comparing the coefficients
of $x^{j-1}$ on both sides of the identity
$$\sum_{k=1}^n(1+x)^{k-1}=\f{(1+x)^n-1}{(1+x)-1},$$
we get a simple proof of the desired identity. \qed

\medskip

\noindent{\it Proof of Theorem 1.1}. Let $\al$ and $\beta$ be the two roots of the equation
$x^2-Ax+B=0$. In view of Lemma 2.1,
$$\align &\sum_{j=1}^{p-1}\f{v_j(A,B)}{jA^j}(-1)^j\bi{p-1}j
\\=&\sum_{j=1}^{p-1}\f{v_j(A,B)}{jA^j}(-1)^j\sum_{k=j}^{p-1}\bi{k-1}{j-1}
\\=&\sum_{k=1}^{p-1}\sum_{j=1}^k\bi{k-1}{j-1}\f{(-1)^jv_j(A,B)}{jA^j}
\\=&\sum_{k=1}^{p-1}\f1k\sum_{j=1}^k\bi kj\(\l(-\f{\al}A\r)^j+\l(-\f{\beta}A\r)^j\)
\\=&\sum_{k=1}^{p-1}\f{(1-\al/A)^k+(1-\beta/A)^k-2}k=\sum_{k=1}^{p-1}\f{\beta^k+\al^k}{kA^k}-2\sum_{k=1}^{p-1}\f1k
\\\eq&\sum_{k=1}^{p-1}\f{v_k(A,B)}{kA^k}\ (\mo\ p^2).
\endalign$$
Since
$$(-1)^k\bi{p-1}k-1\eq-p H_k\ (\mo\ p^2)\ \ \ \t{for}\ k=1,\ldots,p-1,$$
(1.1) follows from the above.

Similarly,
$$\align &(\al-\beta)\sum_{j=1}^{p-1}\f{u_j(A,B)}{jA^j}(-1)^j\bi{p-1}j
\\=&\sum_{j=1}^{p-1}\f{(\al-\beta)u_j(A,B)}{jA^j}(-1)^j\sum_{k=j}^{p-1}\bi{k-1}{j-1}
\\=&\sum_{k=1}^{p-1}\sum_{j=1}^k\bi{k-1}{j-1}\f{(-1)^j(\al-\beta)u_j(A,B)}{jA^j}
\\=&\sum_{k=1}^{p-1}\f1k\sum_{j=1}^k\bi kj\(\l(-\f{\al}A\r)^j-\l(-\f{\beta}A\r)^j\)
\\=&\sum_{k=1}^{p-1}\f{(1-\al/A)^k-(1-\beta/A)^k}k=\sum_{k=1}^{p-1}\f{\beta^k-\al^k}{kA^k}
\\=&(\beta-\al)\sum_{k=1}^{p-1}\f{u_k(A,B)}{kA^k}.
\endalign$$
Thus, if $\Delta=A^2-4B\not=0$ then
$$\sum_{k=1}^{p-1}\f{u_k(A,B)}{kA^k}(1-pH_k)\eq-\sum_{k=1}^{p-1}\f{u_k(A,B)}{kA^k}\ (\mo\ p^2)$$
and hence (1.2) follows.

Now suppose that $\Delta=0$. By induction, $u_k=k(A/2)^{k-1}$ for $k=0,1,2,\ldots$.
Thus
$$\align &\sum_{j=1}^{p-1}\f{u_j(A,B)}{jA^j}(-1)^j\bi{p-1}j
\\=&\sum_{j=1}^{p-1}\f{(-1)^j}{2^{j-1}A}\bi{p-1}j=\f2A\sum_{j=1}^{p-1}\bi{p-1}j\l(-\f12\r)^j
\\=&\f2A\(\l(1-\f12\r)^{p-1}-1\)=\f2A\cdot\f{1-2^{p-1}}{2^{p-1}}
\\=&-\f2A\sum_{j=0}^{p-2}\f{2^j}{2^{p-1}}=-\f 2A\sum_{k=1}^{p-1}\f1{2^k}=-\sum_{k=1}^{p-1}\f{u_k(A,B)}{kA^k}.
\endalign$$
This yields (1.2) with the help of (2.1).
\medskip

So far we have completed the proof of Theorem 1.1. \qed

\proclaim{Lemma 2.2} Let $A,B\in\Z$ and $\Delta=A^2-4B$.
Suppose that $p$ is an odd prime with $p\nmid B\Delta$. Then we have the congruence
$$\l(\f{A\pm\sqrt{\Delta}}2\r)^{p-(\f{\Delta}p)}\eq B^{(1-(\f{\Delta}p))/2}\ \ (\mo\ p)\tag2.2$$
in the ring of algebraic integers.
\endproclaim
\Proof. Both $\al=(A+\sqrt{\Delta})/2$ and $\beta=(A-\sqrt{\Delta})/2$ are roots of the equation
$x^2-Ax+B=0$. Observe that
$$2\al^p\eq 2^p\al^p=(A+\sqrt{\Delta})^p\eq A^p+(\sqrt{\Delta})^p\eq A+\l(\f{\Delta}p\r)\sqrt{\Delta}
\ (\mo\ p).$$
Similarly,
$$2\beta^p\eq A-\l(\f{\Delta}p\r)\sqrt{\Delta}
\ \ (\mo\ p).$$
Thus, if $(\f{\Delta}p)=1$, then
$$\al^{p-1}B=\al^p\beta\eq\al\beta=B\ (\mo\ p)
\ \t{and}\ \beta^{p-1}B=\al\beta^p\eq\al\beta=B\ (\mo\ p),$$
hence $\al^{p-1}\eq1\eq\beta^{p-1}\ (\mo\ p)$.
When $(\f{\Delta}p)=-1$, by the above we have
$$\al^{p+1}=\al\al^p\eq\al\beta=B\ (\mo\ p)\ \ \t{and}\ \ \beta^{p+1}=\beta^p\beta\eq\al\beta=B\ (\mo\ p).$$
This concludes the proof. \qed

\medskip
\noindent {\it Proof of Theorem 1.2}. (i) The equation $x^2-Ax+B=0$ has two roots
$$\al=\f{A+\sqrt{\Delta}}2\ \ \t{and}\ \ \beta=\f{A-\sqrt{\Delta}}2,$$
where $\Delta=A^2-4B$. Also,
$$H_{p-1-k}=H_{p-1}-\sum_{0<j\ls k}\f1{p-j}\eq H_k\ (\mo\ p)\quad\t{for}\ k=0,1,\ldots,p-1.$$
As $(\f{\Delta}p)=1$, with the help of Lemma 2.2, for any $n\in\N$ we have
$$\align \sum_{k=0}^{p-1}v_k(A,B)H_k^n=&\sum_{k=0}^{p-1}v_{p-1-k}(A,B)H_{p-1-k}^n
\\\eq&\sum_{k=0}^{p-1}\l(\al^{p-1-k}+\beta^{p-1-k}\r)H_k^n
\\\eq&\sum_{k=0}^{p-1}\l(\f{\beta^k}{B^k}+\f{\al^k}{B^k}\r)H_k^n
=\sum_{k=0}^{p-1}B^{-k}v_k(A,B)H_k^n\ (\mo\ p).
\endalign$$
This proves (1.5).
Similarly,
$$\align &(\al-\beta)\sum_{k=0}^{p-1}u_k(A,B)H_k^n
=\sum_{k=0}^{p-1}(\al^{p-1-k}-\beta^{p-1-k})H_{p-1-k}^n
\\\eq&\sum_{k=0}^{p-1}\l(\f{\beta^k}{B^k}-\f{\al^k}{B^k}\r)H_k^n
=(\beta-\al)\sum_{k=0}^{p-1}B^{-k}u_k(A,B)H_k^n\ (\mo\ p).
\endalign$$
As $(\al-\beta)^2=\Delta\not\eq0\ (\mo\ p)$, (1.4) follows.

(ii) If $p\eq1\ (\mo\ 3)$ (i.e., $(\f{-3}p)=1$), then by putting $A=\pm1$ and $B=1$ in (1.4) we get
(1.6) and (1.7). Now we prove (1.6) and (1.7) in the case $p\eq2\ (\mo\ 3)$. Observe that
$$\align&\sum_{k=0}^{p-1}(-1)^k\l(\f k3\r)\l(\bi{p-1}k(-1)^k-1\r)
\\=&\sum_{k=0}^{p-1}\bi{p-1}k\l(\f k3\r)-\sum_{k=0}^{p-1}(-1)^k\l(\f k3\r)
\\=&\sum_{k=0}^{p-1}\bi{p-1}k\f{\omega^k-\bar\omega^k}{\sqrt{-3}}
-\sum_{k=0}^{p-1}(-1)^k\f{\omega^k-\bar\omega^k}{\sqrt{-3}}
\\=&\f1{\sqrt{-3}}\l((1+\omega)^{p-1}-(1+\bar\omega)^{p-1}\r)
-\f1{\sqrt{-3}}\l(\f{1+\omega^p}{1+\omega}-\f{1+\bar\omega^p}{1+\bar \omega}\r)
\\=&\f1{\sqrt{-3}}\l((-\omega^2)^{p-1}-(-\omega)^{p-1}-\f{-\omega^{2p}}{-\omega^2}+\f{-\omega^p}{-\omega}\r)=0.
\endalign$$
Also,
$$\align&\sum_{k=0}^{p-1}\l(\f k3\r)\l(\bi{p-1}k(-1)^k-1\r)
\\=&\sum_{k=0}^{p-1}\bi{p-1}k(-1)^k\f{\omega^k-\bar\omega^k}{\sqrt{-3}}-\sum_{j=1}^{(p-2)/3}\sum_{d=0}^2\l(\f{3j-d}3\r)
-\l(\f{p-1}3\r)
\\=&\f1{\sqrt{-3}}\l((1-\omega)^{p-1}-(1-\omega^2)^{p-1}\r)-\l(\f{p-1}3\r)
\\=&\f1{\sqrt{-3}}(1-\omega)^{p-1}\l(1-(-\omega^2)^{p-1}\r)-1
\endalign$$
and hence
$$\align&\sum_{k=0}^{p-1}\l(\f k3\r)\l(\bi{p-1}k(-1)^k-1\r)
\\=&\f{1-\omega^2}{\sqrt{-3}}(1-\omega)^{p-1}-1=\f{-\omega^2}{\sqrt{-3}}(1-\omega)^p-1
\\=&\f{-\omega^2}{\sqrt{-3}}(\sqrt{-3}\ \omega^2)^p-1=-(-3)^{(p-1)/2}\omega^{2+2p}-1
\\=&-\l((-3)^{(p-1)/2}-\l(\f{-3}p\r)\r)
\\\eq&-\f{(\f{-3}p)}2\l((-3)^{p-1}-\l(\f{-3}p\r)^2\r)=\f{3^{p-1}-1}2\ (\mo\ p^2).
\endalign$$
Combining these with (2.1) we immediately obtain (1.6) and (1.7).

Next we show (1.8).  Observe that
$$\align &\sum_{k=0}^{p-(\f p3)}(-1)^k\l(\f k3\r)k
=\sum_{q=1}^{(p-(\f p3))/6}\sum_{r=0}^5(-1)^{6q-r}\l(\f{6q-r}3\r)(6q-r)
\\=&\sum_{q=1}^{(p-(\f p3))/6}((6q-1)+(6q-2)-(6q-4)-(6q-5))=p-\l(\f p3\r)
\endalign$$
and hence
$$\sum_{k=0}^{p-1}(-1)^k\l(\f k3\r)k=\f{1+(\f p3)}2p-\l(\f p3\r).$$
Also,
$$\align&\sum_{k=0}^{p-1}\bi{p-1}k\l(\f k3\r)k=(p-1)\sum_{k=1}^{p-1}\bi{p-2}{k-1}\l(\f k3\r)
\\=&(p-1)\sum_{j=0}^{p-2}\bi{p-2}j\f{\omega^{j+1}-\bar\omega^{j+1}}{\sqrt{-3}}
=\f{p-1}{\sqrt{-3}}\l(\omega(1+\omega)^{p-2}-\bar\omega(1+\bar\omega)^{p-2}\r)
\\=&\f{p-1}{\sqrt{-3}}\l(\omega(-\omega^2)^{p-2}-\omega^2(-\omega)^{p-2}\r)
=\f{p-1}{\sqrt{-3}}\l(\omega^p-\omega^{2p}\r)
=(p-1)\l(\f p3\r).
\endalign$$
Thus
$$\align&\sum_{k=0}^{p-1}(-1)^k\l(\f k3\r)k\l(\bi{p-1}k(-1)^k-1\r)
\\=&(p-1)\l(\f p3\r)-\l(\f{1+(\f p3)}2p-\l(\f p3\r)\r)=\f{(\f p3)-1}2p.
\endalign$$
This implies (1.8) due to (2.1).

Finally we prove (1.9). If $p\eq1\ (\mo\ 4)$ (i.e., $(\f{-4}p)=1$), then (1.4) in the case $A=B=2$ yields (1.9).
Below we assume that $p\eq 3\ (\mo\ 4)$. Note that
$$u_k(2,2)=\f{(1+i)^k-(1-i)^k}{2i}\quad \t{for all}\ k\in\N.$$
Thus
$$\align&2i\sum_{k=0}^{p-1}\bi{p-1}k(-1)^k(2^{-k}+1)u_k(2,2)
\\=&\sum_{k=0}^{p-1}\bi{p-1}k((-2)^{-k}+(-1)^k)\l((1+i)^k-(1-i)^k\r)
\\=&\l(1-\f{1+i}2\r)^{p-1}-\l(1-\f{1-i}2\r)^{p-1}+(1-(1+i))^{p-1}-(1-(1-i))^{p-1}
\\=&\l(\f{1-i}2\r)^{p-1}-\l(\f{1+i}2\r)^{p-1}+(-i)^{p-1}-i^{p-1}
\\=&\l(\f{-2i}4\r)^{(p-1)/2}-\l(\f{2i}4\r)^{(p-1)/2}=2i\f{i^{(p+1)/2}}{2^{(p-1)/2}}
\endalign$$
and hence
$$\sum_{k=0}^{p-1}\bi{p-1}k(-1)^k(2^{-k}+1)u_k(2,2)=\f{(-1)^{(p+1)/4}}{2^{(p-1)/2}}.\tag2.3$$
Also,
$$\align&2i\sum_{k=0}^{p-1}(2^{-k}+1)u_k(2,2)
\\=&\sum_{k=0}^{p-1}(2^{-k}+1)\l((1+i)^k-(1-i)^k\r)
\\=&\f{1-((1+i)/2)^p}{1-(1+i)/2}-\f{1-((1-i)/2)^p}{1-(1-i)/2}+\f{1-(1+i)^p}{1-(1+i)}-\f{1-(1-i)^p}{1-(1-i)}
\\=&(1+i)-2\l(\f{1+i}2\r)^{p+1}-\l((1-i)-2\l(\f{1-i}2\r)^{p+1}\r)
\\&+i-i(1+i)(1+i)^{p-1}+i-i(1-i)(1-i)^{p-1}
\\=&2i-2\l(\f{2i}4\r)^{(p+1)/2}+2\l(\f{-2i}4\r)^{(p+1)/2}
\\&+2i+(1-i)(2i)^{(p-1)/2}-(1+i)(-2i)^{(p-1)/2}
\\=&4i+2(2i)^{(p-1)/2}=2i\l(2+i^{(p-3)/2}2^{(p-1)/2}\r)
\endalign$$
and hence
$$\sum_{k=0}^{p-1}(2^{-k}+1)u_k(2,2)=2-(-1)^{(p+1)/4}2^{(p-1)/2}.\tag2.4$$
Combining (2.1), (2.3) and (2.4) we obtain
$$\align &-p\sum_{k=0}^{p-1}(2^{-k}+1)u_k(2,2)H_k
\\\eq&\f{(-1)^{(p+1)/4}}{2^{(p-1)/2}}-2+(-1)^{(p+1)/4}2^{(p-1)/2}
\\\eq&\f{(-1)^{(p+1)/4}}{2^{(p-1)/2}}\l(2^{(p-1)/2}-(-1)^{(p+1)/4}\r)^2\eq0\ (\mo\ p^2)
\endalign$$
since
$$\l(\f 2p\r)=(-1)^{(p^2-1)/8}=(-1)^{(p+1)/4\times(p-1)/2}=(-1)^{(p+1)/4}.$$
Therefore (1.9) holds.

\smallskip

 By the above, we have completed the proof of Theorem 1.2. \qed

\heading{3. Proof of Theorem 1.3}\endheading

\proclaim{Lemma 3.1} Let $A,B\in\Z$. Let $p$ be an odd prime with $(\f Bp)=1$.
Suppose that $b^2\eq B\ (\mo\ p)$ where $b\in\Z$.
 Then
$$u_{(p-1)/2}(A,B)\eq\cases0\ (\mo\ p)&\t{if}\ (\f{A^2-4B}p)=1,\\\f1b(\f{A-2b}p)\ (\mo\ p)&\t{if}\ (\f{A^2-4B}p)=-1;
\endcases$$
and
$$u_{(p+1)/2}(A,B)\eq\cases(\f{A-2b}p)\ (\mo\ p)&\t{if}\ (\f{A^2-4B}p)=1,
\\0\ (\mo\ p)&\t{if}\ (\f{A^2-4B}p)=-1.
\endcases$$
\endproclaim
\Proof. The congruences are known results, see, e.g., [S]. \qed

\proclaim{Lemma 3.2} Let $u_n=u_n(1,4)$ for $n\in\N$. Then, for any prime $p>5$ we have
$$u_p-2^{p-1}\l(\f p{15}\r)\eq 2^{(\f p{15})-2}u_{p-(\f p{15})}\ (\mo\ p^2).\tag3.1$$
\endproclaim
\Proof. The two roots
$$\al=\f{1+\sqrt{-15}}2\quad\ \t{and}\ \quad\beta=\f{1-\sqrt{-15}}2$$
of the equation $x^2-x+4=0$ are algebraic integers.
Clearly
$$-15u_p=(\al-\beta)^2u_p=(\al-\beta)(\al^p-\beta^p)\eq(\al-\beta)^{p+1}=(-15)^{(p+1)/2}\ (\mo\ p)$$
and hence
$$u_p\eq (-15)^{(p-1)/2}\eq\l(\f{-15}p\r)=\l(\f p{15}\r)\ (\mo\ p).$$
Also,
$$v_p=\al^p+\beta^p\eq(\al+\beta)^p=1\ (\mo\ p),$$
where  $v_n$ refers to $v_n(1,4)$.
(In fact, both $u_p(A,B)$ and $v_p(A,B)$ modulo $p$ are known for any $A,B\in\Z$.)
By induction, $u_n+v_n=2u_{n+1}$ for any $n\in\N$.

{\it Case} 1. $(\f p{15})=1.$ In this case,
$$4u_{p-1}=u_p-u_{p+1}=\f{u_p-v_p}2\eq\f{1-1}2=0\ (\mo\ p)$$
and
$$v_{p-1}=2u_p-u_{p-1}\eq 2\eq 2^p\ (\mo\ p).$$
Since
$$\align (v_{p-1}-2^p)(v_{p-1}+2^p)=&(\al^{p-1}+\beta^{p-1})^2-4(\al\beta)^{p-1}
\\=&(\al^{p-1}-\beta^{p-1})^2=-15u_{p-1}^2\eq0\ (\mo\ p^2),
\endalign$$
we must have
$v_{p-1}\eq 2^p\ (\mo\ p^2)$ and hence
$$2u_p=u_{p-1}+v_{p-1}\eq u_{p-1}+2^p\ (\mo\ p^2).$$

{\it Case} 2. $(\f p{15})=-1$. In this case,
$$2u_{p+1}=u_p+v_p\eq -1+1=0\ (\mo\ p)$$
and
$$v_{p+1}=2u_{p+2}-u_{p+1}=u_{p+1}-8u_p\eq8\eq 2^{p+2}\ (\mo\ p).$$
As
$$\align (v_{p+1}-2^{p+2})(v_{p+1}+2^{p+2})=&(\al^{p+1}+\beta^{p+1})^2-4(\al\beta)^{p+1}
\\=&(\al^{p+1}-\beta^{p+1})^2=-15u_{p+1}^2\eq0\ (\mo\ p^2),
\endalign$$
we must have
$v_{p+1}\eq 2^{p+2}\ (\mo\ p^2)$ and hence
$$\align 8u_p=&2(u_{p+1}-u_{p+2})=2u_{p+1}-(u_{p+1}+v_{p+1})
\\=&u_{p+1}-v_{p+1}\eq u_{p+1}-2^{p+2}\  (\mo\ p^2).
\endalign$$

Combining the above, we immediately obtain the desired result. \qed

\medskip
\noindent{\it Proof of Theorem 1.3}.
Set $\da=(1-(\f p{15}))/2$. Then
$$\sum_{k=0}^{p-1}\f{u_{k+\delta}}{2^k}H_k\eq
\sum_{k=0}^{p-1}\f{u_{k+\delta}}{2^k}\cdot\f{1-(-1)^k\bi{p-1}k}p\ (\mo\ p).$$
So it suffices to show
$$\sum_{k=0}^{p-1}u_{k+\delta}2^{p-1-k}\eq \sum_{k=0}^{p-1}\bi{p-1}ku_{k+\delta}(-2)^{p-1-k}\ (\mo\ p^2),\tag3.2$$
which implies (1.11) and (1.12) in the cases $\da=0,1$ respectively.
Recall that
$$u_{k+\da}=\f{\al^{k+\da}-\beta^{k+\da}}{\al-\beta},$$
where
$$\al=\f{1+\sqrt{-15}}2\quad\t{and}\quad\beta=\f{1-\sqrt{-15}}2$$
are the two roots of the equation $x^2-x+4=0$.
Since
$$\sum_{k=0}^{p-1}x^ky^{p-1-k}=\f{x^p-y^p}{x-y}\quad\t{and}\quad\sum_{k=0}^{p-1}\bi{p-1}kx^ky^{p-1-k}=(x+y)^{p-1},$$
(3.2) can be rewritten as follows:
$$\aligned&\f1{\al-\beta}\(\al^{\delta}\f{\al^{p}-2^{p}}{\al-2}-\beta^\da\f{\beta^{p}-2^{p}}{\beta-2}\)
\\\eq&\f{\al^{\da}(\al-2)^{p-1}-\beta^\da(\beta-2)^{p-1}}{\al-\beta}\ (\mo\ p^2).
\endaligned\tag3.3$$
Note that $(\al-2)(\beta-2)=4+\al\beta-2(\al+\beta)=4+4-2=6$ and
$$\align&\al^\da(\al^{p}-2^{p})(\beta-2)-\beta^\da(\al-2)(\beta^{p}-2^{p})
\\=&(2^{p}-\al^{p})(\al^{\da+1}+\al^\da)+(\beta^{\da+1}+\beta^\da)(\beta^{p}-2^{p})
\\=&2^{p}(\al^{\da+1}-\beta^{\da+1}+\al^\da-\beta^\da)-(\al^{p+\da+1}-\beta^{p+\da+1})-(\al^{p+\da}-\beta^{\beta+\da})
\\=&(\al-\beta)\l(2^{p}(u_{\da+1}+u_\da)-(u_{p+\da+1}+u_{p+\da})\r)
\\=&(\al-\beta)\l(2^{p+\da}-2u_{p+\da}+4u_{p+\da-1}\r).
\endalign$$
So the left-hand side of (3.3) coincides with
$$\f{2^{p+\da}-2u_{p+\da}+4u_{p+\da-1}}6=\f{2^{p+\da-1}-u_{p+\da}+2u_{p+\da-1}}3.$$

Since $(\al-2)^2=-3\al$ and $(\beta-2)^2=-3\beta$, we have
$$\align &\f{\al^\da(\al-2)^{p-1}-\beta^\da(\beta-2)^{p-1}}{\al-\beta}
\\=&\f{\al^\da(-3\al)^{(p-1)/2}-\beta^\da(-3\beta)^{(p-1)/2}}{\al-\beta}
\\=&(-3)^{(p-1)/2}\f{\al^{(p-1)/2+\da}-\beta^{(p-1)/2+\da}}{\al-\beta}=(-3)^{(p-1)/2}u_{(p-(\f p{15}))/2}.
\endalign$$
Applying Lemma 3.1 with $A=1$ and $B=4$, we get that
$$u_{(p-(\f p{15}))/2}\eq0\ (\mo\ p)\ \ \t{and}\ \ u_{(p+(\f p{15}))/2}\eq\l(\f{-3}p\r)2^{((\f{p}{15})-1)/2}\ (\mo\ p).
\tag3.4$$

For any $n\in\N$ we have
$$u_{2n}=\f{\al^{2n}-\beta^{2n}}{\al-\beta}=\f{\al^n-\beta^n}{\al-\beta}(\al^n+\beta^n)=u_nv_n.$$
If $(\f p{15})=1$, then by (3.4) we have $u_{(p-1)/2}\eq0\ (\mo\ p)$ and
$$v_{(p-1)/2}=2u_{(p+1)/2}-u_{(p-1)/2}\eq2\l(\f{-3}p\r)\ (\mo\ p),$$
hence
$$\align u_{p-1}=&u_{(p-1)/2}v_{(p-1)/2}
\\\eq&2\l(\f{-3}p\r)u_{(p-1)/2}\eq 2(-3)^{(p-1)/2}u_{(p-1)/2}\ (\mo\ p^2).
\endalign$$
If $(\f p{15})=-1$, then by (3.4) we have $u_{(p+1)/2}\eq0\ (\mo\ p)$ and
$$\align v_{(p+1)/2}=&2u_{(p+3)/2}-u_{(p+1)/2}=u_{(p+1)/2}-8u_{(p-1)/2}
\\\eq&-8\l(\f{-3}p\r)2^{-1}=-4\l(\f{-3}p\r)\ (\mo\ p),
\endalign$$
hence
$$\align u_{p+1}=&u_{(p+1)/2}v_{(p+1)/2}
\\\eq&-4\l(\f{-3}p\r)u_{(p+1)/2}\eq -4(-3)^{(p-1)/2}u_{(p+1)/2}\ (\mo\ p^2).
\endalign$$
Thus the right-hand side of (3.3) is congruent to $u_{p-(\f p{15})}/(2(-2)^\da)$ mod $p^2$.

By the above, (3.3) is equivalent to the following congruence
$$\f{2^{p+\da-1}-u_{p+\da}+2u_{p+\da-1}}3\eq\f{u_{p-(\f p{15})}}{2(-2)^\da}\ (\mo\ p^2).\tag3.5$$
If $(\f p{15})=1$, then $\da=0$, and hence (3.5) reduces to the congruence
$$2(2^{p-1}-u_p+2u_{p-1})\eq 3u_{p-1}\pmod {p^2}$$
which is equivalent to (3.1) since $(\f p{15})=1$.
When $(\f p{15})=-1$, we have $\da=1$ and hence (3.5) can be rewritten as
$$-4(2^p-u_{p+1}+2u_p)\eq 3u_{p+1}\pmod{p^2}$$
which follows from (3.1) since $(\f p{15})=-1$.
This concludes the proof.  \qed

\Ack. The author is grateful to the two referees for their helpful comments.


\medskip

 \widestnumber\key{BPQ}

 \Refs

\ref\key BPQ\by A. T. Benjamin, G. O. Preston and J. J. Quinn\paper A Stirling encounter with harmonic numbers
\jour Math. Mag. \vol 75\yr 2002\pages 95--103\endref

\ref\key BB\by  D. Borwein and J. M. Borwein\paper On an intriguing integral and some series related to $\zeta(4)$
\jour Proc. Amer. Math. Soc.\vol 123\yr 1995\pages 1191--1198\endref

\ref\key C\by S. W. Coffman\paper Problem 1240 and Solution: An infinite series with harmonic numbers
\jour Math. Mag.\vol 60\yr 1987\pages 118--119\endref

\ref\key GKP\by R. L. Graham, D. E. Knuth and O. Patashnik
 \book Concrete Mathematics\publ 2nd ed., Addison-Wesley, New York\yr 1994\endref

\ref\key HT\by C. Helou and G. Terjanian\paper On Wolstenholme's theorem and its converse
\jour J. Number Theory \vol 128\yr 2008\pages 475--499\endref

\ref\key R\by P. Ribenboim\book The Book of Prime Number Records
\publ Springer, New York, 1989\endref

\ref\key S\by Z. H. Sun\paper Values of Lucas sequences modulo primes
\jour Rocky Mount. J. Math. \vol 33\yr 2003\pages 1123--1145\endref


\ref\key Su\by Z. W. Sun\paper Arithmetic theory of harmonic numbers
\jour Proc. Amer. Math. Soc.\pages in press\endref


\ref\key W\by J. Wolstenholme\paper On certain properties of prime numbers
\jour Quart. J. Math.\vol 5\yr 1862\pages 35--39\endref

\endRefs
\enddocument